\documentclass{article}

\usepackage[applemac]{inputenc} 
\usepackage[dvips]{graphicx}
\usepackage{amsmath} 
\usepackage{amssymb} 

\usepackage{colortbl}
\usepackage{soul} 

\usepackage{pstcol}
\usepackage{pst-node}
\usepackage{pstricks}
\usepackage{fancybox}
\usepackage{calc}
\usepackage{truncate}
\usepackage{comment}

\makeatletter
\def\@maketitle{%
  \vbox to 2.6in{%
    \hsize\textwidth
    \linewidth\hsize
    \vspace*{1.5cm}
    \centering
    {\bfseries\huge \@title \par}
    \vskip 2em
    {\large \begin{tabular}[t]{c}\@author \end{tabular}\par}
    \vfill}    \vspace*{1.0cm}
}
\renewcommand\section{\@startsection {section}{1}{\z@}%
     {.7\baselineskip plus\baselineskip}{.5\baselineskip}
                                   {\normalfont\Large\bfseries}}
\renewcommand\section{\@startsection {section}{1}{\z@}%
      {.5\baselineskip\@plus.7\baselineskip}{.3\baselineskip}%
                                   {\normalfont\Large\bfseries}}
\renewcommand\subsection{\@startsection{subsection}{2}{\z@}%
       {.5\baselineskip\@plus.7\baselineskip}{.3\baselineskip}%
                                   {\normalfont\large\bfseries}}
\renewcommand\subsubsection{\@startsection{subsubsection}{3}{\z@}%
      {.5\baselineskip\@plus.7\baselineskip}{.3\baselineskip}%
                                     {\normalfont\normalsize\bfseries}}
\makeatother

\renewenvironment{abstract}%
  {\normalfont
    \list{}{\labelwidth0pt
      \leftmargin0pt \rightmargin\leftmargin
      \listparindent\parindent \itemindent0pt
      \parsep0pt
      
    }%
    \item[\hskip\labelsep\bfseries\abstractname\enspace --] \itshape%
}{%
  \endlist}

\newcommand{\keywordsname}{Keywords}
\newenvironment{keywords}%
  {\normalfont
    \list{}{\labelwidth0pt
      \leftmargin0pt \rightmargin\leftmargin
      \listparindent\parindent \itemindent0pt
      \parsep0pt
      }%
    \item[\hskip\labelsep\bfseries\keywordsname:]}{\endlist}

\begin{document}
\pagestyle{plain} 

\title{
Partial ordering of hyper-powersets\\
 and matrix representation \\
of belief functions within DSmT}

\author{Jean Dezert\\
ONERA, Dtim/ied\\
29 Av. de la  Division Leclerc \\
92320 Ch\^{a}tillon, France.\\
Jean.Dezert@onera.fr\\
\and
Florentin Smarandache\\
Department of Mathematics\\
University of New Mexico\\
Gallup, NM 87301, U.S.A.\\
smarand@unm.edu}
\date{}

\maketitle

\begin{abstract}
In this paper, we examine several issues for ordering or partially ordering elements of hyper-powertsets  involved in the recent theory of plausible, uncertain and paradoxical reasoning (DSmT) developed by the authors. We will show the benefit of some of these issues to obtain a nice and useful matrix representation of belief functions.
\end{abstract}


\section{Introduction}

The Dezert-Smarandache theory (DSmT for short) of plausible, uncertain and paradoxical reasoning \cite{Dezert_2002,Dezert_2002b,Dezert_2003,Smarandache_2002} is a generalization of the classical Dempster-Shafer theory (DST) \cite{Shafer_1976} which allows to formally combine any types of sources of information (rational, uncertain or paradoxical). The DSmT is able to solve complex data/information fusion problems where the DST usually fails, specially when conflicts (paradoxes) between sources become large and when the refinement of the frame of discernment 
$\Theta$ is inaccessible because of the vague, relative and imprecise nature of elements of $\Theta$ (see \cite{Dezert_2003} for justification and examples). The foundation of DSmT is based on the definition of the hyper-powerset $D^\Theta$ (or free distributive lattice on $n$ generators) of a general frame of discernment $\Theta$.  $\Theta$ must be considered as a set $\{\theta_{1},\ldots,\theta_{n}\}$ of $n$  elements considered as exhaustive which cannot be 
precisely defined and separated, so that no refinement of $\Theta$ into a new larger set $\Theta_{ref}$ of disjoint elementary hypotheses is possible in contrast with the classical Shafer's model on which is based the DST. We have already presented in a companion paper \cite{Dezert_2003d}, how to easily generate all elements of $D^\Theta$ using the property of isotone Boolean functions. In this paper, we focus our attention, on how to order them in a clever way in order to get a very interesting matrix representation of belief functions defined over $D^\Theta$.The DSmT deals directly with paradoxical/conflicting sources of information into this new formalism and proposes a new and very simple (associative and commutative) rule of combination for conflicting sources of informations (corpus/bodies of evidence). Some interesting results based on DSmT approach can be found in \cite{Tchamova_2003,Corgne_2003}.
Before going deeper into the DSmT it is necessary to briefly present first the foundations of the DST and DSmT for a better understanding of the important differences between these two theories based on Shafer model and DSm model.

\section{Short presentation of the DST}

\subsection{The Shafer's model}

The Shafer's model assumes that  the frame of discernment of the problem under consideration is a set $\Theta=\{\theta_1,\theta_2,\ldots,\theta_n\}$ of $n$ {\it{exhaustive}} and {\it{exclusive elementary}} hypothesis $\theta_i$. Such model implicitly imposes that an ultimate refinement of the problem is always possible so that $\theta_i$ can be well precisely defined/identified in such a way  that we are sure that they are exclusive and exhaustive. From this model,  Shafer defines a basic belief assignment (bba) $m (.): 2^\Theta \rightarrow  [0, 1]$  associated to a given body of evidence $\mathcal{B}$ by adding the following constraints to $m(.)$
\begin{equation}
m(\emptyset)=0  \qquad \text{and}\qquad     \sum_{A\in 2^\Theta} m(A) = 1                            
\end{equation}
\noindent
where $2^\Theta$ is called the {\it{powerset}} of $\Theta$, i.e. the set of all subsets of $\Theta$. From any bba, one then defines the belief and plausibility functions for all  $A\subseteq\Theta$ as
\begin{equation}
\text{Bel}(A) = \sum_{B\in 2^\Theta, B\subseteq A} m(B)
\label{Belg}
\end{equation}
\begin{equation}
\text{Pl}(A) = \sum_{B\in 2^\Theta, B\cap A\neq\emptyset} m(B)=1- \text{Bel}(\bar{A})
\label{Plg}
\end{equation}

\subsection{The Dempster's combination rule}

Let $\text{Bel}_1(.)$ and $\text{Bel}_2(.)$ be two belief functions over the same frame of discernment $\Theta$ and their corresponding bba $m_1(.)$ and $m_2(.)$ provided by two {\it{distinct}} bodies of evidence $\mathcal{B}_1$ and $\mathcal{B}_2$. Then the combined global belief function $\text{Bel}(.)= \text{Bel}_1(.)\oplus \text{Bel}_2(.)$ is obtained by combining the information granules $m_1(.)$ and $m_2(.)$ through the  Dempster's rule of combination $[m_{1}\oplus m_{2}](\emptyset)=0$ and $\forall B\neq\emptyset \in 2^\Theta$ as follows
 \begin{equation}
[m_{1}\oplus m_{2}](B) = 
\frac{\sum_{X\cap Y=B}m_{1}(X)m_{2}(Y)}{1-\sum_{X\cap Y=\emptyset} m_{1}(X) m_{2}(Y)} 
\label{eq:DSR}
 \end{equation}
 
The notation $\sum_{X\cap Y=B}$ represents the sum over all $X, Y \in 2^\Theta$ such that $X\cap Y=B$. 
The orthogonal sum $m (.)\triangleq [m_{1}\oplus m_{2}](.)$ is considered as a basic belief assignment if and only if the denominator in equation \eqref{eq:DSR} is non-zero. The term $k_{12}\triangleq \sum_{X\cap Y=\emptyset} m_{1}(X) m_{2}(Y)$ is called degree of conflict between the sources $\mathcal{B}_1$ and $\mathcal{B}_2$. When $k_{12}=1$,  the orthogonal sum $m (.)$ does not exist and the bodies of evidences $\mathcal{B}_1$ and $\mathcal{B}_2$ are said to be in {\it{full contradiction}}. Such a case can arise when there exists $A \subset \Theta$ such that $\text{Bel}_1(A) =1$ and $\text{Bel}_2(\bar{A}) = 1$. Same kind of trouble can occur also with the {\it{Optimal Bayesian Fusion Rule}} (OBFR) \cite{Dezert_2001a,Dezert_2001b}. 

\subsection{Alternatives to Dempter's rule}

The Shafer's model and the DST is attractive for the {\it{Data Fusion community}} because it gives a nice mathematical model for ignorance and it includes the Bayesian theory as a special case \cite{Shafer_1976} (p.4). Although very appealing, the DST presents nevertheless some important weaknesses and limitations because of its model itself, the theoretical justification of the Dempster's rule of combination but also because of our confidence to trust the result of Dempster's rule of combination specially when the conflict becomes important between sources ($k_{12} \nearrow 1$).The {\it{a posteriori}} justification of the Dempster's rule of combination has been brought by the Smets axiomatic of the Transferable Belief Model (TBM) in \cite{Smets_1994}. 
But recently, we must also emphasize here that an infinite number of possible rules of combinations can be built from the Shafer's model  following ideas initially proposed in \cite{Lefevre_2002} and corrected here as follows:
\begin{itemize}
\item one first has to compute $m(\emptyset)$ by
$$ m(\emptyset) \triangleq \sum_{A\cap B =\emptyset}m_1(A)m_2(B) $$
\item then one redistributes $m(\emptyset)$ on all $A\subseteq \Theta$ with some given positive coefficients $w_m(A)$ such that $\sum_{A\subseteq \Theta} w_m(A)=1$ according to
\begin{equation}
\begin{cases}
w_m(\emptyset)m(\emptyset) \rightarrow m(\emptyset)\\
m(A) + w_m(A)m(\emptyset) \rightarrow m(A), \forall A\neq\emptyset
\end{cases}
\label{eq:CEV}
\end{equation}
\end{itemize}
The particular choice of the set of coefficients $w_m(.)$ provides a particular rule of combination.
Actually there exists an infinite number of possible rules of combination. Some rules can be better justified than others depending on their ability or not to preserve associativity and commutativity properties of the combination. It can be easily shown in  \cite{Lefevre_2002}  that such general procedure provides all existing rules developed in the literature from the Shafer's model as alternative to the primeval Dempster's rule of combination depending on the choice of coefficients $w(A)$. As example the Dempster's rule of combination can be obtained from \eqref{eq:CEV} by choosing $w_m(\emptyset)=0$ and $ w_m(A)=m(A)/(1-m(\emptyset))$ for all $A\neq\emptyset$. The Yager's rule of combination is obtained by choosing $w_m(\Theta)=1$ while the "Smets' rule of combination" is obtained by choosing $w_m(\emptyset)=1$ and thus accepting the possibility to deal with bba such that $m(\emptyset)>0$.

\subsection{Matrix calculus for belief functions}

As rightly emphasized recently by Smets in \cite{Smets_2002}, the mathematic of belief functions is often cumbersome  because of the many summations symbols and all its subscripts involved in equations.  This renders equations very difficult to read and understand at first sight and might discourage potential readers for their complexity. Actually, this is just an appearance because most of the operations encountered in DST with belief functions and basic belief assignments $m(.)$ are just simple linear operations and can be easily represented using matrix notation and be handled by elementary matrix calculus. We just focus here our presentation on the matrix representation of the relationship between a basic belief assignment $m(.)$ and its associated belief function $\text{Bel}(.)$. A nice and more complete presentation of matrix calculus for belief functions can be found in \cite{Kennes_1991,Kennes_1992,Smets_2002}. One important aspect for the simplification of matrix representation and calculus in DST, concerns the choice of the order of the elements of the powerset $2^\Theta$. The order of elements of $2^\Theta$ can be chosen arbitrarily actually, and it can be easily seen by denoting $\mathbf{m}$ the bba vector of size $2^n\times 1$ and $\text{{\bf{Bel}}}$ its corresponding belief vector of same size, that the set of equations \eqref{Belg} holding for all $A\subseteq\Theta$ is strictly equivalent to the following general matrix equation
\begin{equation}
\text{{\bf{Bel}}}= \mathbf{BM}\cdot\mathbf{m} \quad \Leftrightarrow \quad \mathbf{m}= \mathbf{BM}^{-1}\cdot \text{{\bf{Bel}}}
\label{eq:BMDST}
\end{equation}
\noindent
where the internal structure of $\mathbf{BM}$ depends on the choice of the order for enumerating the elements of $2^\Theta$.
But it turns out that the simplest ordering based on the enumeration of integers from 0 to $2^{n}-1$ expressed as $n$-binary strings with the lower bit on the right (LBR) (where $n=\vert \Theta\vert$) to characterize all the elements of powerset, is the most efficient solution and best encoding method for matrix calculus and for developing efficient algorithms in MatLab\footnote{Matlab is a trademark of The MathWorks, Inc.}  or similar programming languages \cite{Smets_2002}. By choosing the basic increasing binary enumeration (called  {\it{bibe system}}), one obtains a very nice recursive algorithm on the dimension $n$ of $\Theta$ for computing the matrix $\mathbf{BM}$. The computation of $\mathbf{BM}$ for $\vert\Theta\vert=n$ is just obtained  from the iterations up to $i+1=n$ of the recursive relation \cite{Smets_2002} starting with $\mathbf{BM}_{0}\triangleq[1]$ and where $ \mathbf{0}_{i+1}$ denotes the zero-matrix of size $(i+1)\times(i+1)$,
\begin{equation}
\mathbf{BM}_{i+1}=
\begin{bmatrix} 
\mathbf{BM}_{i} & \mathbf{0}_{i+1}\\
\mathbf{BM}_{i} & \mathbf{BM}_{i}
\end{bmatrix}
\end{equation}
\noindent $\mathbf{BM}$ is a binary unimodular matrix ($\text{det}(\mathbf{BM})=\pm 1$).  $\mathbf{BM}$ is moreover triangular inferior and symmetrical with respect to its antidiagonal.\\

\noindent
{\bf{Example for $\Theta=\{\theta_{1},\theta_{2},\theta_{3}\}$}}\\
The {\it{bibe system}} gives us the following order for elements of $2^\Theta=\{\alpha_0,\ldots,\alpha_7\}$:

\begin{center}
\begin{tabular}{ll}
$\alpha_0\equiv 000\equiv \emptyset$  & $\alpha_1\equiv 001 \equiv \theta_{1} $ \\
$\alpha_2\equiv 010 \equiv \theta_{2}$ & $\alpha_3\equiv 011 \equiv  \theta_{1}\cup\theta_{2}$ \\
$\alpha_4\equiv 100\equiv \theta_{3}$  & $\alpha_5\equiv 101 \equiv \theta_{1}\cup\theta_{3}$\\
$\alpha_6\equiv 110 \equiv \theta_{2}\cup\theta_{3}$ & $\alpha_7\equiv 111 \equiv  \theta_{1}\cup\theta_{2}\cup\theta_{3}\equiv  \Theta$
\end{tabular}
\end{center}
\noindent
Each element $\alpha_i$ of $2^\Theta$ is a 3-bits string. With this bibe system, on has  $\mathbf{m}=[m(\alpha_0),\ldots,m(\alpha_7)]'$ and  $\text{{\bf{Bel}}}=[\text{Bel}((\alpha_0),\ldots,\text{Bel}((\alpha_7)]'$. The expressions of the matrix $\mathbf{BM}_3$ and its inverse ${\mathbf{BM}_3}^{-1}$ are given by
\begin{equation*}
\mathbf{BM}_3=\begin{bmatrix} 
\textcolor{red}{1} & 0 & 0 & 0 & 0 & 0 &0 & 0\\
\textcolor{red}{1} & \textcolor{red}{1} & 0 & 0 & 0 & 0 &0 & 0\\
\textcolor{red}{1} & 0 & \textcolor{red}{1} & 0 & 0 & 0 &0 & 0\\
\textcolor{red}{1} & \textcolor{red}{1} & \textcolor{red}{1} & \textcolor{red}{1} & 0 & 0 &0 & 0\\
\textcolor{red}{1} & 0 & 0 & 0 & \textcolor{red}{1} & 0 &0 & 0\\
\textcolor{red}{1} & \textcolor{red}{1} & 0 & 0 & \textcolor{red}{1} & \textcolor{red}{1} &0 & 0\\
\textcolor{red}{1} & 0 & \textcolor{red}{1} & 0 & \textcolor{red}{1} & 0 &\textcolor{red}{1} & 0\\
\textcolor{red}{1} & \textcolor{red}{1} & \textcolor{red}{1} & \textcolor{red}{1} & \textcolor{red}{1} & \textcolor{red}{1} &\textcolor{red}{1} & \textcolor{red}{1}
\end{bmatrix}
\end{equation*}
\begin{equation*}
{\mathbf{BM}_3}^{-1}=
\begin{bmatrix} 
\textcolor{red}{1} & 0 & 0 & 0 & 0 & 0 &0 & 0\\
\textcolor{blue}{-1} & \textcolor{red}{1} & 0 & 0 & 0 & 0 &0 & 0\\
\textcolor{blue}{-1} & 0 & \textcolor{red}{1} & 0 & 0 & 0 &0 & 0\\
\textcolor{red}{1} & \textcolor{blue}{-1} & \textcolor{blue}{-1} & \textcolor{red}{1} & 0 & 0 &0 & 0\\
\textcolor{blue}{-1} & 0 & 0 & 0 & \textcolor{red}{1} & 0 &0 & 0\\
\textcolor{red}{1} & \textcolor{blue}{-1} & 0 & 0 & \textcolor{blue}{-1} & \textcolor{red}{1} &0 & 0\\
\textcolor{red}{1} & 0 & \textcolor{blue}{-1} & 0 & \textcolor{blue}{-1} & 0 &\textcolor{red}{1} & 0\\
\textcolor{blue}{-1} & \textcolor{red}{1} & \textcolor{red}{1} & \textcolor{blue}{-1} & \textcolor{red}{1} & \textcolor{blue}{-1} &\textcolor{blue}{-1} & \textcolor{red}{1}
\end{bmatrix}
\end{equation*}

\section{A short DSmT presentation}

\subsection{The DSm model}

The development of the Dezert-Smarandache Theory (DSmT) of plausible, uncertain, and paradoxical reasoning comes from the necessity to overcome, for a wide class of problems, the two following inherent limitations of the DST which are closely related with the acceptance of the third middle excluded principle, i.e.
\begin{enumerate}
\item[(C1)] - the DST considers a discrete and finite frame of discernment 
$\Theta$ based on a set of exhaustive and exclusive elementary elements $\theta_i$.
\item[(C2)] - the bodies of evidence are assumed independent and provide their own belief function on the powerset $2^\Theta$ but with {\it{same interpretation}} for $\Theta$. 
\end{enumerate}

The relaxation of constraints (C1) and (C2) seems necessary for a wide class of fusion problems due to the possible vague, imprecise and paradoxical nature of the elements of $\Theta$. By accepting the third middle, we can easily handle the possibility to deal directly with a new kinds of elements with respect to those belonging to the Shafer's model. This is the DSm model. A wider class of interesting fusion problems  can then be attacked by the DSmT. The relaxation of the constraint (C1) can be justified since, in many problems (see example in \cite{Dezert_2003}), the elements of $\Theta$ generally correspond only to imprecise/vague notions and concepts so that no refinement of $\Theta$ satisfying the first constraint is actually possible. The relaxation of (C2) is also justified since, in general, the {\sl{same}} frame $\Theta$ may be interpreted differently by the distinct sources of evidence. Some subjectivity on the information provided by a source is almost unavoidable. In most of cases, the sources of evidence provide their beliefs about some hypotheses only with respect to their own worlds of knowledge, experiences, feelings, senses without reference to the (inaccessible) absolute truth of the space of possibilities and without any probabilistic background argumentations. The DSmT  includes the possibility to deal with evidences arising from different sources of information which don't have access to absolute interpretation of the elements $\Theta$ under consideration. The DSmT can be interpreted as a general and direct extension of Bayesian theory and the Dempster-Shafer theory in the 
following sense. Let $\Theta=\{\theta_{1},\theta_{2}\}$ be the simplest frame of 
discernment involving only two elementary hypotheses (with no  additional assumptions on $\theta_{1}$ and $\theta_{2}$), then 
\begin{itemize}
\item the probability theory deals with basic probability assignments (bpa)
$m(.)\in [0,1]$ such that $$m(\theta_{1})+m(\theta_{2})=1$$
\item the DST deals with bba $m(.)\in [0,1]$ such that 
$$m(\theta_{1})+m(\theta_{2})+m(\theta_{1}\cup\theta_{2})=1$$
\item the DSmT theory deals with
generalized bba $m(.)\in [0,1]$ such that 
$$m(\theta_{1})+m(\theta_{2})+m(\theta_{1}\cup\theta_{2})+m(\theta_{1}\cap\theta_{2})=1$$
\end{itemize}

\subsection{DSm model versus Shafer's model}

The Shafer's model considers that the frame $\Theta$ of the problem under consideration is a set of finite exhaustive and exclusive elements $\theta_i$ and requires in some way a refinement in order to choose/select $\theta_i$ as exclusive. The DSm model can be viewed as the model opposite to the Shafer's model where none of the $\theta_i$ are considered exclusive. This DSm model is justified in a wide class of fusion problems when the intrinsic nature of the elements of $\Theta$ to be manipulated is such that $\Theta$ is not refinable at all into exclusive and precise subsets \cite{Dezert_2003}. The DSmT can then deal with elements/concepts which have possibly (but not necessary) continuous and/or relative interpretation to  the corpus of evidences like, by example, the relative notions of smallness/tallness, beauty/ugliness, pleasure/pain, heat/coldness, even the notion of colors (due to the continuous spectrum of the light),  etc. None of these notions or concepts can be clearly refined/separated in an absolute manner so that they cannot be considered as exclusive and we cannot also define precisely what their conjunctions are. Their interpretations/estimations through the bba mechanism given by any corpus of evidence is always built from its own (limited) knowledge/experience and senses. Between these two extreme models, there exists a finite number of {\it{DSm-hybrid}} models for which some integrity constraints (by forcing some potential conjunctions to be impossible, i.e. equal to the empty set) between some elements of $\theta$ can be introduced depending on the hybrid-nature of the problem. The DSm model can then be viewed as the most free model and the Shafer's model as the most restrictive one.  The DSmT has been developed up to now only for the DSm model but application of the DSmT for DSm-hybrid models is under investigation.

\subsection{Notion of hyper-powerset $D^\Theta$}

One of the cornerstones of the DSmT is  the notion of hyper-powerset which is defined as follows.
Let $\Theta=\{\theta_{1},\ldots,\theta_{n}\}$ be a set of $n$ 
elements which cannot be precisely defined and separated so that no 
refinement of $\Theta$ in a new larger set $\Theta_{ref}$ of disjoint elementary 
hypotheses is possible (we abandon here the Shafer's model). The {\sl{hyper-powerset}} $D^\Theta$ is defined as the set of all composite propositions built from elements of $\Theta$ with $\cup$ and $\cap$ ($\Theta$ generates $D^\Theta$ under operators $\cup$ and $\cap$)
operators such that 
\begin{enumerate}
\item $\emptyset, \theta_1,\ldots, \theta_n \in D^\Theta$.
\item  If $A,B \in D^\Theta$, then $A\cap B\in D^\Theta$ and $A\cup B\in D^\Theta$.
\item No other elements belong to $D^\Theta$, except those obtained by using rules 1 or 2.
\end{enumerate}
The dual (obtained by switching $\cup$ and $\cap$ in expressions) of $D^\Theta$ is itself.  There are elements in $D^\Theta$ which are self-dual (dual to themselves), for example $\alpha_8$ for the case when $n=3$ in the example below.
The cardinality of $D^\Theta$ is majored by 
$2^{2^n}$ when $\text{Card}(\Theta)=\vert\Theta\vert=n$. The generation 
of hyper-powerset $D^\Theta$ is closely related with the famous Dedekind's problem on enumerating the 
set of isotone Boolean functions \cite{Dezert_2003d}. The cardinality of $D^\Theta$ for $n=\vert\Theta\vert=0,1,2,3, ...$ follows the sequence of Dedekind's numbers 1,2,5,19,167,7580,... \cite{Dezert_2003d}. From a general frame of discernment $\Theta$, we define a map $m(.): 
D^\Theta \rightarrow [0,1]$ associated to a given body of evidence $\mathcal{B}$ 
which can support paradoxical information, as follows
\begin{equation*}
m(\emptyset)=0 \qquad \text{and}\qquad \sum_{A\in D^\Theta} m(A) = 1 
\end{equation*}
The quantity $m(A)$ is called $A$'s {\sl{generalized basic belief 
assignment}} (gbba) or the generalized basic belief mass for $A$.
The belief and plausibility functions are defined in almost the same manner as within the DST, i.e.
\begin{equation}
\text{Bel}(A) = \sum_{B\in D^\Theta, B\subseteq A} m(B)
\label{BelgDSmT}
\end{equation}
\begin{equation}
\text{Pl}(A) = \sum_{B\in D^\Theta, B\cap A\neq\emptyset} m(B)
\end{equation}
These definitions are compatible with the DST definitions when the sources of 
information become uncertain but rational (they do not support paradoxical 
information). We still have $\forall A\in D^\Theta, \text{Bel}(A)\leq \text{Pl}(A)$.

\subsection{The DSm rule of combination}

The  DSm rule of combination $m(.)\triangleq [m_{1}\oplus m_{2}](.)$ of two distinct (but potentially paradoxical) sources of evidences $\mathcal{B}_{1}$ and  $\mathcal{B}_{2}$ over the same 
general frame of discernment $\Theta$ with belief functions $\text{Bel}_{1}(.)$ and 
 $\text{Bel}_{2}(.)$ associated with general information granules $m_{1}(.)$ and $m_{2}(.)$ is  given by $\forall C\in D^\Theta$,
 \begin{equation*}
m(C) = 
 \sum_{A,B\in D^\Theta, A\cap B=C}m_{1}(A)m_{2}(B)
 \label{JDZT}
 \end{equation*}
Since $D^\Theta$ is closed under $\cup$ and $\cap$ operators, this new rule 
of combination guarantees that $m(.): D^\Theta \rightarrow [0,1]$ is a proper general information granule. This rule of combination is commutative and associative 
and can always be used for the fusion of paradoxical or rational sources of 
information (bodies of evidence). It is important to note that any fusion of sources of information 
generates either uncertainties, paradoxes or {\sl{more generally both}}. This is intrinsic to the 
general fusion process itself. The theoretical justification of the DSm rule can be found in \cite{Dezert_2003}. A network representation of this DSm rule of combination can be found in \cite{Dezert_2003d}.\\

\section{Ordering elements of hyper-powerset for matrix calculus}

As within the DST framework,  the order of the elements of $D^\Theta$ can be  arbitrarily chosen.  We denote the Dedekind number or order $n$ as $d(n)\triangleq \vert D^\Theta\vert$ for $n=\vert\Theta\vert$. We denote also $\mathbf{m}$ the gbba vector of size $d(n)\times 1$ and $\text{{\bf{Bel}}}$ its corresponding belief vector of the same size. The  set of equations \eqref{BelgDSmT} holding for all $A\in D^\Theta$ is then strictly equivalent to the following general matrix equation
\begin{equation}
\text{{\bf{Bel}}}= \mathbf{BM}\cdot\mathbf{m} \quad \Leftrightarrow \quad \mathbf{m}= \mathbf{BM}^{-1}\cdot \text{{\bf{Bel}}}
\label{eq:BMDSmT}
\end{equation}

Note the similarity between these relations with the previous ones \eqref{eq:BMDST}. The only difference resides in the size of vectors 
$\text{{\bf{Bel}}}$ and $ \mathbf{m}$ and the size of matrix $\mathbf{BM}$ and their components. We explore in the following sections the possible choices for ordering (or partially ordering) the elements of hyper-powerset $D^\Theta$, to obtain an interesting matrix structure of $\mathbf{BM}$ matrix. Only three issues are examined and briefly presented in the sequel. The first method is based on the direct enumeration of elements of $D^\Theta$ according to their recursive generation via the algorithm for generating all isotone Boolean functions presented in \cite{Dezert_2003d}. The second (partial) ordering method is based on the notion of DSm cardinality which will be introduced in section 4.2. The last and most interesting solution proposed for partial ordering over $D^\Theta$ is obtained by introducing the notion of intrinsic informational strength $s(.)$ associated to each element of hyper-powerset.

\subsection{Order based on the enumeration of isotone Boolean functions}

We have presented in \cite{Dezert_2003d} a recursive algorithm based on isotone Boolean functions for generating $D^\Theta$. Here is briefly the principle of the method. Let consider $\Theta=\{\theta_1,\ldots,\theta_n\}$ satisfying the DSm model and the Dezert-Smarandache order $\mathbf{u}_n$ of the Smarandache's codification of parts of Venn diagram $\Theta$ with $n$ partially overlapping elements $\theta_i$, $i=1,\ldots, n$ (see \cite{Dezert_2003d} for details about Smarandache's codification). All the elements $\alpha_i$ of $D^\Theta$ can then be obtained by the very simple linear equation \cite{Dezert_2003d}
\begin{equation}
\mathbf{d}_n= \mathbf{D}_n\cdot\mathbf{u}_n
\end{equation}
\noindent where $\mathbf{d}_n\equiv[\alpha_0\equiv\emptyset,\alpha_1,\ldots,\alpha_{d(n)-1}]'$ is the vector of elements of $D^\Theta$, $\mathbf{u}_n$ is the proper codification vector and $D_n$ a particular binary matrix. The final result $\mathbf{d}_n$ is obtained from the previous {\it{matrix product}} after identifying $(+,\cdot)$ with $(\cup,\cap)$ operators, $0\cdot x$ with $\emptyset$ and $1\cdot x$ with $x$). $D_n$ is actually a binary matrix corresponding to isotone (i.e. non-decreasing) Boolean functions obtained by applying recursively the steps
(starting with $\mathbf{D}_0^c=[0 \, 1]'$)
\begin{itemize}
\item $\mathbf{D}_n^c$ is built from $\mathbf{D}_{n-1}^c$ by adjoining to each row $\mathbf{r}_i$ of $\mathbf{D}_{n-1}^c$ any row $\mathbf{r}_j$ of $\mathbf{D}_{n-1}^c$ such that $\mathbf{r}_i \cup \mathbf{r}_j=\mathbf{r}_j$. Then $\mathbf{D}_n$ is obtained by removing the first column and the last line of $\mathbf{D}_n^c$.
\end{itemize}

\noindent
{\bf{Example for $\Theta=\{\theta_{1},\theta_{2},\theta_{3}\}$}}
\begin{equation*}
\underbrace{
\begin{bmatrix}
    \alpha_0\\
    \alpha_1\\
    \alpha_2\\
    \alpha_3\\
    \alpha_4\\
    \alpha_5\\
    \alpha_6\\
    \alpha_7\\
    \alpha_8\\
    \alpha_9\\
    \alpha_{10}\\
    \alpha_{11}\\
    \alpha_{12}\\
    \alpha_{13}\\
    \alpha_{14}\\
    \alpha_{15}\\
    \alpha_{16}\\
    \alpha_{17}\\
    \alpha_{18}
\end{bmatrix}
}_{\mathbf{d}_3}
=
\underbrace{
\begin{bmatrix}
    0  &   0    & 0   &  0   &  0   &  0   &  0\\
     0 &    0   &  0  &   0  &   0  &   0  &   1\\
     0  &   0   &  0   &  0  &   0  &   1   &  1\\
     0  &   0  &   0  &   0  &   1  &   0   &  1\\
     0  &   0  &   0  &   0  &   1  &   1   &  1\\
     0   &  0  &   0  &   1  &   1  &   1   &  1\\
     0  &   0   &  1   &  0  &   0  &   0   &  1\\
     0  &   0  &   1  &   0  &   0  &   1   &  1\\
     0  &   0  &   1  &   0  &   1  &   0  &   1\\
     0  &   0   &  1  &   0  &   1  &   1  &   1\\
     0  &   0   &  1  &   1   &  1  &   1  &   1\\
     0  &   1  &   1 &    0  &   0   &  1  &   1\\
     0  &   1  &   1  &   0  &   1   &  1  &   1\\
     0  &   1  &   1  &   1   &  1  &   1  &   1\\
     1  &   0  &   1  &   0  &   1  &   0 &    1\\
     1  &   0  &   1  &   0  &   1  &   1   &  1\\
     1  &   0  &   1  &   1  &   1  &   1 &    1\\
     1   &  1  &   1  &   0   &  1  &   1  &   1\\
     1  &   1  &   1  &   1  &   1  &   1  &   1
\end{bmatrix}
}_{\mathbf{D}_3}
\cdot
\underbrace{
\begin{bmatrix}
   <1>\\
    <2>\\
    <12>\\
    <3>\\
    <13>\\
     <23>\\
     <123>\\
\end{bmatrix}
}_{\mathbf{u}_3}
 \end{equation*}

Hence we finally get (after simple algebraic simplifications) the following irreducible elements for $D^\Theta$
\begin{equation*}
\begin{array}{l}
\alpha_i   \; \text{(from the isotone Boolean functions alg.)}  \\
\hline
\alpha_0\triangleq\emptyset                                                 \\
\alpha_1\triangleq\theta_1\cap\theta_2\cap\theta_3       \\
\alpha_2\triangleq\theta_2\cap\theta_3                               \\
\alpha_3\triangleq\theta_1\cap\theta_3                                \\
\alpha_4\triangleq(\theta_1\cup\theta_2)\cap\theta_3       \\
\alpha_5\triangleq\theta_3                                \\
\alpha_6\triangleq\theta_1\cap\theta_2                               \\
\alpha_7\triangleq(\theta_1\cup\theta_3)\cap\theta_2       \\
\alpha_8\triangleq(\theta_2\cup\theta_3)\cap\theta_1       \\
\alpha_9\triangleq[(\theta_1\cap\theta_2)\cup\theta_3] \cap(\theta_1\cup\theta_2)  \\
\alpha_{10}\triangleq(\theta_1\cap\theta_2)\cup\theta_3      \\
\alpha_{11}\triangleq\theta_2                               \\
\alpha_{12}\triangleq(\theta_1\cap\theta_3)\cup\theta_2      \\
\alpha_{13}\triangleq(\theta_2\cup\theta_3)                             \\
\alpha_{14}\triangleq\theta_1                               \\
\alpha_{15}\triangleq(\theta_2\cap\theta_3)\cup\theta_1      \\
\alpha_{16}\triangleq(\theta_1\cup\theta_3)                             \\
\alpha_{17}\triangleq(\theta_1\cup\theta_2)                             \\
\alpha_{18}\triangleq(\theta_1\cup\theta_2\cup\theta_3)       
\end{array}
\end{equation*}

We denote $r^{iso}(\alpha_i)$ the position of $\alpha_i$ into the column vector $\mathbf{d}_n$ obtained from the previous enumeration/generation system. Such system  provides a total order over $D^\Theta$ defined $\forall \alpha_i,\alpha_j \in D^\Theta$ as $\alpha_i\prec\alpha_j$ ($\alpha_i$ precedes $\alpha_j$) if and only if $r^{iso}(\alpha_i) < r^{iso}(\alpha_j)$. 
Based on this order, the $\mathbf{BM}$ matrix involved in \eqref{eq:BMDSmT} presents unfortunately no particular interesting structure.
We have thus to look  for better solutions for ordering the elements of hyper-powersets.

\subsection{Ordering with the DSm cardinality}

A second possibility for ordering the elements of $D^\Theta$ is to (partially) order them by their increasing {\it{DSm cardinality}}. The {\it{DSm cardinality}} of any element $A \in D^\Theta$, denoted $\mathcal{C}_\mathcal{M}(A)$,
corresponds to the number of parts of $A$ in the Venn diagram of the problem (model $\mathcal{M}$) taking into account the set of integrity constraints (if any), i.e. all the possible intersections due to the nature of the elements $\theta_i$. This {\it{intrinsic cardinality}} depends on the model $\mathcal{M}$.  $\mathcal{M}$ is the model that contains $A$ which depends on the
dimension of Venn diagram, (i.e. the number of sets $n=\vert \Theta \vert$ under consideration),
and on the number of non-empty intersections in this diagram. One has $1 \leq \mathcal{C}_\mathcal{M}(A) \leq 2^n-1$. $\mathcal{C}_\mathcal{M}(A)$ must not be confused with the classical cardinality $\vert A \vert$ of a given set $A$ (i.e. the number of its distinct elements) - that's why a new notation is necessary here. \\

In the (general) case of the free-model $\mathcal{M}^f$ (i.e. the DSm model) where all conjunctions are non-empty, one has for intersections:
\begin{itemize}
\item $\mathcal{C}_{\mathcal{M}^f}(\theta_1)= \ldots=\mathcal{C}_{\mathcal{M}^f}(\theta_n)=2^{n-1}$
\item $\mathcal{C}_{\mathcal{M}^f}(\theta_i\cap\theta_j)= 2^{n-2}$ for $n\geq 2$
\item $\mathcal{C}_{\mathcal{M}^f}(\theta_i\cap\theta_j\cap\theta_k)= 2^{n-3}$ for $n\geq 3$
\end{itemize}
It can be proved by induction that for $1\leq m \leq n$, one has $\mathcal{C}_{\mathcal{M}^f}(\theta_{i_1}\cap\theta_{i_2}\cap\ldots\cap\theta_{i_m})= 2^{n-m}$. For the cases $n=1,2,3,4$, this formula can be checked on the corresponding Venn diagrams. Let's consider this formula true for $n$ sets, and prove it for $n+1$ sets (when all intersections/conjunctions are considered non-empty).  From the Venn diagram of $n$ sets, we can get a Venn diagram with $n+1$ sets if one draws a closed curve that cuts each of the $2^{n}-1$ parts of the previous diagram (and, as a consequence, divides each part into two disjoint subparts).  Therefore, the number of parts of each intersection is doubling when passing from a diagram of dimension $n$ to a diagram of dimension $n+1$.  Q.e.d.\\

In the case of the free-model $\mathcal{M}^f$, one has for unions:
\begin{itemize}
\item $\mathcal{C}_{\mathcal{M}^f}(\theta_i\cup\theta_j)= 3(2^{n-2})$ for $n\geq 2$
\item $\mathcal{C}_{\mathcal{M}^f}(\theta_i\cup\theta_j\cup\theta_k)= 7(2^{n-3})$ for $n\geq 3$
\end{itemize}
It can be proved also by induction that for $1\leq m \leq n$, one has $\mathcal{C}_{\mathcal{M}^f}(\theta_{i_1}\cup\theta_{i_2}\cup\ldots\cup\theta_{i_m})= (2^m-1)(2^{n-m})$. The proof is similar to the previous one, and keeping in mind that passing from a Venn diagram of dimension $n$ to a dimension $n+1$, all each part that forms the union $\theta_i\cap\theta_j\cap\theta_k$ will be split into two disjoint parts, hence the number of parts is doubling.\\

For other elements $A$ in $D^\Theta$, formed by unions and intersections, the close-form for $\mathcal{C}_{\mathcal{M}^f}(A)$ seems more complicated to obtain. But from the generation algorithm of $D^\Theta$ (see \cite{Dezert_2003d} for details), DSm cardinal of a set $A$ from $D^\Theta$ is exactly equal to the sum of its coefficients in the $\mathbf{u}_n$ basis, i.e. the sum of its row elements in the $\mathbf{D}_n$ matrix, which is actually very easy to compute by programming. The {\it{DSm cardinality}} plays in important role in the definition of the Generalized Pignistic Transform (GPT) for the construction of subjective/pignistic probabilities of elements of $D^\Theta$ for decision-making \cite{Dezert_Smarandache_Daniel_2003}.\\

If one imposes a constraint that a set $B$ from $D^\Theta$ is empty, then one suppresses the columns corresponding to the parts which compose $B$ in the $\mathbf{D}_n$ matrix and the row of $B$ and the rows of all elements  of $D^\Theta$ which are subsets of $B$, getting a new matrix ${\mathbf{D}'}_n$ which represents a new model $\mathcal{M}'$.
In the $\mathbf{u}_n$ basis, one similarly suppresses the parts that form $B$, and now this basis has the dimension $2^n-1-\mathcal{C}_{\mathcal{M}}(B)$.\\

\noindent {\it{Example with $\mathcal{M}^f$}}: Consider the 3D case $\Theta=\{\theta_1,\theta_2,\theta_3\}$ with the free-model $\mathcal{M}^f$  corresponding to the following Venn diagram (where $<i>$ denotes the part 
which belongs to $\theta_i$ only, $<ij>$ denotes the part which belongs to $\theta_i$ and $\theta_j$ only, etc; this is the Smarandache's codification \cite{Dezert_2003d}).

\begin{figure}[hbtp]
\centering
\includegraphics[width=4cm]{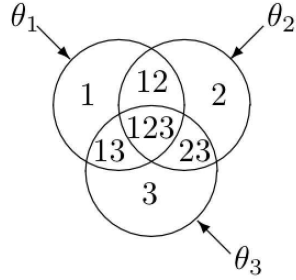}
\caption{Venn Diagram for $\mathcal{M}^f$ }
 \label{fig:2}
 \end{figure}

The corresponding {\it{partial ordering}} for elements of $D^\Theta$ is then summarized in the following table:

\begin{equation*}
\begin{array}{lcl}
A\in D^\Theta                     & \mathcal{C}_{\mathcal{M}^f}(A) \\
\hline
\emptyset                                                                   & 0 \\
\theta_1\cap\theta_2\cap\theta_3       & 1 \\
\theta_1\cap\theta_2                              & 2 \\
\theta_1\cap\theta_3                              & 2 \\
\theta_2\cap\theta_3                              & 2 \\
(\theta_1\cup\theta_2)\cap\theta_3     & 3 \\
(\theta_1\cup\theta_3)\cap\theta_2      & 3 \\
(\theta_2\cup\theta_3)\cap\theta_1      & 3 \\
\theta_1                                                      & 4 \\
\theta_2                                                      & 4 \\
\theta_3                                                      & 4 \\
\{(\theta_1\cap\theta_2)\cup\theta_3\} \cap(\theta_1\cup\theta_2) & 4 \\
(\theta_1\cap\theta_2)\cup\theta_3       &  5 \\
(\theta_1\cap\theta_3)\cup\theta_2       & 5 \\
(\theta_2\cap\theta_3)\cup\theta_1       & 5 \\
\theta_1\cup\theta_2                             & 6 \\
\theta_1\cup\theta_3                             & 6 \\
\theta_2\cup\theta_3                             & 6 \\
\theta_1\cup\theta_2\cup\theta_3         & 7  \\
\end{array}
\end{equation*}

Note that this partial ordering differs from the one described in the previous section and doesn't properly catch the intrinsic informational structure/strength of elements since by example $\{(\theta_1\cap\theta_2)\cup\theta_3\} \cap(\theta_1\cup\theta_2)$ and $\theta_1$ have the same DSm cardinal although they don't look similar because the part $<1>$ in $\theta_1$ belongs only to $\theta_1$ but none of the parts of $\{(\theta_1\cap\theta_2)\cup\theta_3\} \cap(\theta_1\cup\theta_2)$ belongs to only one part of some $\theta_i$. A better ordering function is then necessary to catch the intrinsic informational structure of elements of $D^\Theta$. This is the purpose of the next section.\\

\noindent {\it{Example with another model}}: Consider now the same 3D case with the model $\mathcal{M}\neq\mathcal{M}^f$ in which we force all possible conjunctions to be empty, but $\theta_1\cap\theta_2$ according to the following Venn diagram.

\begin{figure}[hbtp]
\centering
\includegraphics[width=4cm]{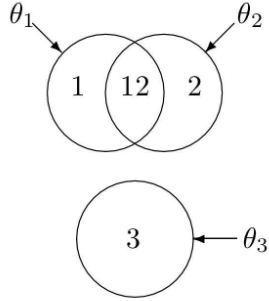}
\caption{Venn Diagram for $\mathcal{M}$ }
 \label{fig:2}
 \end{figure}

The corresponding {\it{partial ordering}} for elements of $D^\Theta$, taking into account the constraints of this model, is then summarized in the following table:

\begin{equation*}
\begin{array}{lcl}
A\in D^\Theta                     & \mathcal{C}_{\mathcal{M}}(A) \\
\hline
\emptyset                                                 & 0 \\
\theta_1\cap\theta_2                              & 1 \\
\theta_3                                                      & 1 \\
\theta_1                                                      & 2 \\
\theta_2                                                      & 2 \\
\theta_1\cup\theta_2                             & 3 \\
\theta_1\cup\theta_3                             & 3 \\
\theta_2\cup\theta_3                             & 3 \\
\theta_1\cup\theta_2\cup\theta_3       & 4  \\
\end{array}
\end{equation*}

The partial ordering of $D^\Theta$ based on DSm cardinality does not provide in general an efficient solution to get an interesting structure for the $\mathbf{BM}$ matrix involved in \eqref{eq:BMDSmT}, contrarily to the structure obtained by Smets in the DST framework (sec. 2.4). The partial ordering presented in the sequel will however allow us to get such nice structure for the matrix calculus of belief functions.

\subsection{Ordering based on the intrinsic informational content}

As already reported, the DSm cardinality is insufficient to catch the intrinsic informational content of each element $d_i$ of $D^\Theta$. A better approach to obtain this, is based on the following new function $s(.)$, which describes  the {\sl{intrinsic information strength}} of any  $d_i\in D^\Theta$. A previous, but cumbersome, definition of $s(.)$ had been proposed in our previous works \cite{Dezert_2002b,Dezert_2003} but  it was difficult to handle and questionable with respect to the formal  equivalent (dual) representation of elements belonging to $D^\Theta$. We propose here a new solution for $s(.)$, based on a very simple and natural geometrical interpretation of the relationships between the parts of the Venn diagram belonging to each $d_i\in D^\Theta$. All the values of the $s(.)$ function (stored into a vector $\mathbf{s}$) over $D^\Theta$ are defined by the following equation:
\begin{equation}
\mathbf{s}= \mathbf{D}_n\cdot \mathbf{w}_n
\end{equation}
\noindent with $\mathbf{s}\triangleq [s(d_0)\: \ldots \:s(d_p)]'$ where $p$ is the cardinal of $D^\Theta$ for the model $\mathcal{M}$ under consideration. $p$ is equal to the Dedekind's number $d(n)-1$ if the free-model $\mathcal{M}^f$ is chosen for $\Theta=\{\theta_1,\ldots,\theta_n\}$. $\mathbf{D}_n$ is the {\it{hyper-powerset generating matrix}}. The components $w_i$ of vector $\mathbf{w}_n$ are obtained from the components of the {\it{Dezert-Smarandache encoding basis}} vector
$\mathbf{u}_n$ as follows  (see \cite{Dezert_2003d} for definitions and details about $\mathbf{D}_n$ and $\mathbf{u}_n$) :
\begin{equation}
w_i\triangleq 1/l(u_i)
\end{equation}
\noindent where $l(u_i)$ is the length of Smarandache's codification $u_i$ of the part of the Venn diagram of the model $\mathcal{M}$, i.e the number of symbols involved in the codification. For example, if $u_i=<123>$, then $l(u_i)=3$ just because only three symbols 1, 2, and 3 enter in the codification $u_i$, thus $w_i = 1/3$.\\

From this new DSm ordering function s(.) we can partially order all the elements  $d_i\in D^\Theta$ by the increasing values of $s(.)$.\\

\noindent{\it{Example for $\Theta=\{\theta_1,\theta_2\}$ with the free-model $\mathcal{M}^f$}}:\\

In this simple case, the DSm ordering of $D^\Theta$ is given by
\begin{equation*}
\begin{array}{l|l}
\alpha_i                         & s(\alpha_i) \\
\hline
\alpha_0=\emptyset                          &s(\alpha_0)=0\\
\alpha_1=\theta_1\cap\theta_2      & s(\alpha_1)=1/2\\
\alpha_2=\theta_1                             & s(\alpha_2)=1+1/2\\
\alpha_3=\theta_2                             & s(\alpha_3)=1+1/2\\
\alpha_4=\theta_1\cup\theta_2    &s(\alpha_{4})=1+1+1/2
\end{array}
\end{equation*}

Based on this ordering, it can be easily verified that the matrix calculus of the beliefs  $\text{{\bf{Bel}}}$ from $\mathbf{m}$ by equation \eqref{eq:BMDSmT}, is equivalent to
$$
\underbrace{
\begin{bmatrix} 
\text{Bel}(\emptyset)\\
\text{Bel}(\theta_1\cap\theta_2)\\
\text{Bel}(\theta_1)\\
\text{Bel}(\theta_2)\\
\text{Bel}(\theta_1\cup\theta_2)
\end{bmatrix}
}_{\text{{\bf{Bel}}}}
=
\underbrace{
\begin{bmatrix} 
\textcolor{blue}{1} & \textcolor{blue}{0} & \textcolor{blue}{0} & \textcolor{blue}{0} & \textcolor{blue}{0} \\
\textcolor{blue}{1} &  \textcolor{red}{1} & 0 & 0 &  0\\
\textcolor{blue}{1} &  \textcolor{red}{1} &  \textcolor{red}{1} & 0 & 0 \\
\textcolor{blue}{1} & \textcolor{red}{1} & 0 & \textcolor{red}{1} & 0 \\
\textcolor{blue}{1} & \textcolor{red}{1} & \textcolor{red}{1} & \textcolor{red}{1} & \textcolor{red}{1} 
\end{bmatrix}
}_{\mathbf{BM}_2}
\underbrace{
\begin{bmatrix} 
m(\emptyset)\\
m(\theta_1\cap\theta_2)\\
m(\theta_1)\\
m(\theta_2)\\
m(\theta_1\cup\theta_2)
\end{bmatrix}
}_{\mathbf{m}}
$$

\noindent where the $ \mathbf{BM}_2$ matrix has a interesting structure (triangular inferior and unimodular properties, $det(\mathbf{BM}_2)=det(\mathbf{BM}_2^{-1})=1$). Conversely,  the calculus of the generalized  basic belief assignment $\mathbf{m}$ from beliefs  $\text{{\bf{Bel}}}$ will be obtained by the inversion of the previous linear system of equations
$$
\underbrace{
\begin{bmatrix} 
m(\emptyset)\\
m(\theta_1\cap\theta_2)\\
m(\theta_1)\\
m(\theta_2)\\
m(\theta_1\cup\theta_2)
\end{bmatrix}
}_{\mathbf{m}}
=
\underbrace{
\begin{bmatrix} 
\textcolor{blue}{1} & \textcolor{blue}{0} & \textcolor{blue}{0} & \textcolor{blue}{0} & \textcolor{blue}{0} \\
\textcolor{blue}{-1} &  \textcolor{red}{1} & 0 & 0 &  0\\
\textcolor{blue}{0} &  \textcolor{red}{-1} &  \textcolor{red}{1} & 0 & 0 \\
\textcolor{blue}{0} & \textcolor{red}{-1} & 0 & \textcolor{red}{1} & 0 \\
\textcolor{blue}{0} & \textcolor{red}{1} & \textcolor{red}{-1} & \textcolor{red}{-1} & \textcolor{red}{1} 
\end{bmatrix}
}_{\mathbf{MB}_2=\mathbf{BM}_2^{-1}}
\underbrace{
\begin{bmatrix} 
\text{Bel}(\emptyset)\\
\text{Bel}(\theta_1\cap\theta_2)\\
\text{Bel}(\theta_1)\\
\text{Bel}(\theta_2)\\
\text{Bel}(\theta_1\cup\theta_2)
\end{bmatrix}
}_{\text{{\bf{Bel}}}}
$$

\noindent{\it{Example for $\Theta=\{\theta_1,\theta_2,\theta_3\}$ with the free-model $\mathcal{M}^f$}}:\\

In this more complicated case, the DSm ordering of $D^\Theta$ is now given by

\begin{figure}[hbtp]
\centering
\includegraphics[width=8cm]{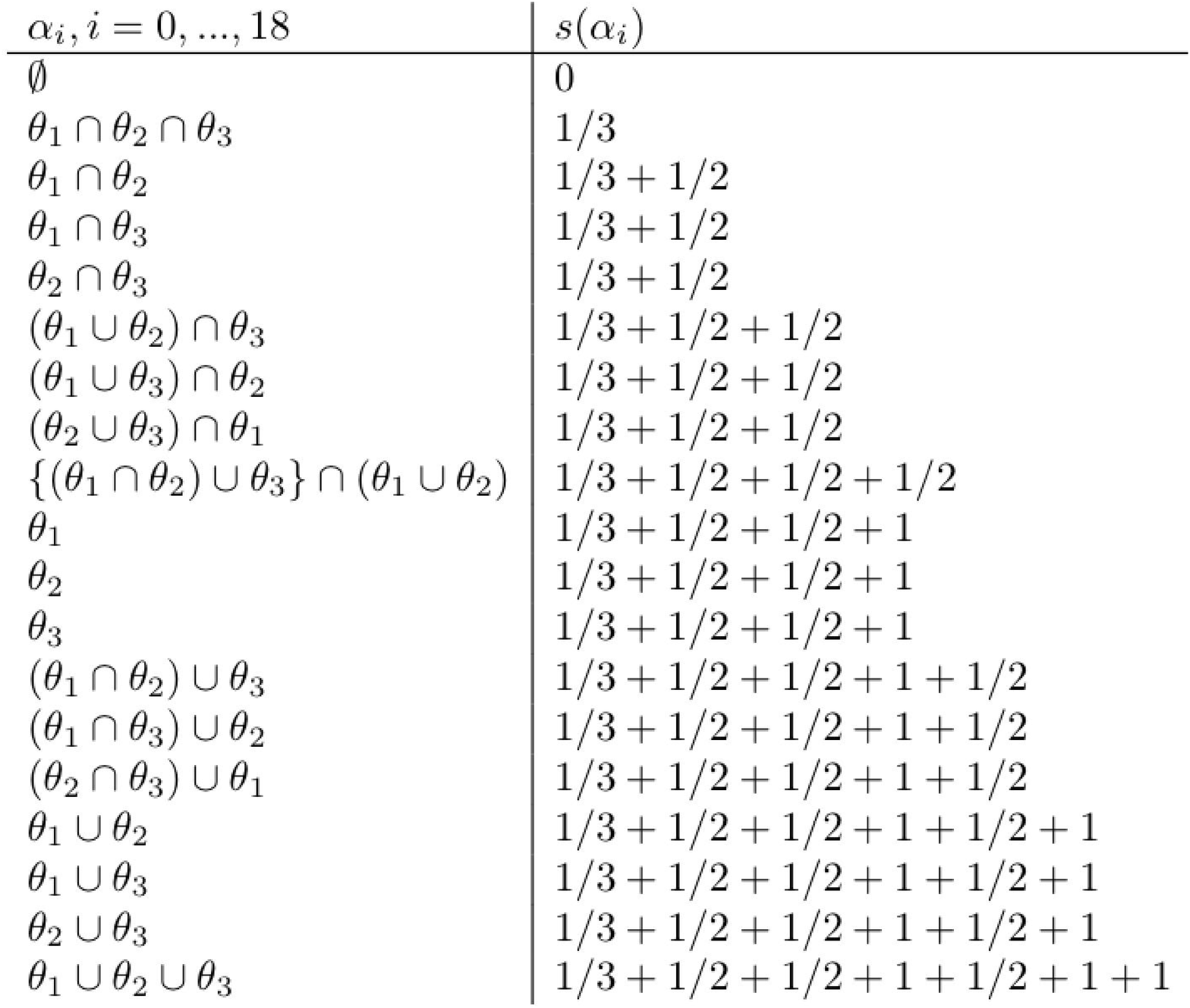}
 \end{figure}

The structure of the matrix $\mathbf{BM}_3$ associated to this ordering is given by
\begin{figure}[hbtp]
\centering
\includegraphics[width=8cm]{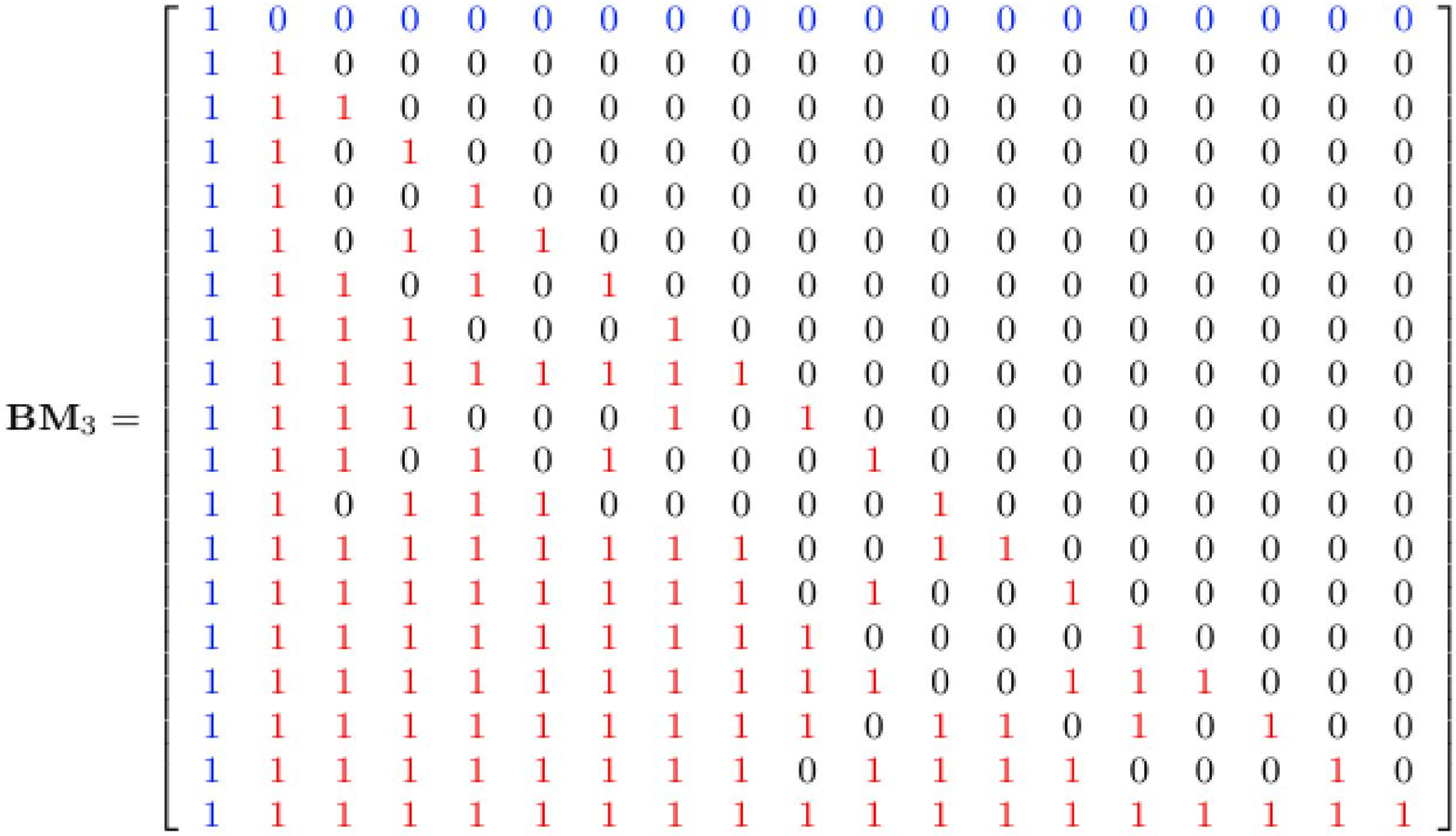}
 \end{figure}

The order for elements generating  the same value of $s(.)$ can be chosen arbitrarily and doesn't change the structure of the matrix $\mathbf{BM}_3$. That's why only a partial order is possible from $s(.)$. It can be verified that $\mathbf{BM}_3$ holds also the same previous interesting matrix structure properties and that $det(\mathbf{BM}_3)=det(\mathbf{BM}_3^{-1})=1$. Similar structure can be shown for problems of higher dimensions ($n>3$).\\

Although a nice structure for matrix calculus of belief functions has been obtained in this work, and conversely to the recursive construction of $\mathbf{BM}_n$ in DST framework, a recursive algorithm (on dimension $n$) for the construction of $\mathbf{BM}_n$ from $\mathbf{BM}_{n-1}$ has not yet be found and is still an open problem for further research.

\section{Conclusion}

A recent theory of plausible, uncertain, and paradoxical reasoning (DSmT) has been developed by the authors to deal with conflicting/paradoxist sources of information which could not be solved by Dempster-Shafer theory of evidence (DST).  DSm rule of combining works for any kind of sources of information (certain, uncertain, paradoxist) depending on each particular model (problem), whereas DS rule of combining fails when the degree of conflict is high.  In order to obtain an easy matrix representation of the belief functions in the DSmT, we need to better order the elements of hyper-powerset $D^\Theta$, that's why we propose in this paper three such orderings: first, using the direct enumeration of isotone Boolean functions, second, based on the DSm cardinality, and third, and maybe the most interesting, by introducing the intrinsic informational strength function $s(.)$ constructed in the DSm encoding basis.


\begin{thebibliography}{99}

\bibitem{Corgne_2003}
Corgne S., Hubert-Moy L., Dezert J., Mercier G., \emph{Land cover change prediction with a new theory of plausible and paradoxical reasoning}, Proc. of the 6th Int. Conf. on inf. fusion (Fusion 2003), Cairns, Australia, july 8-11, 2003.

\bibitem{Dezert_2001a}
Dezert J.,  \emph{Optimal Bayesian fusion of multiple unreliable classifiers}, Proc. of the 4th Int. Conf. on inf. fusion (Fusion 2001), Montréal, Canada, Aug. 8-11, 2001.

\bibitem{Dezert_2001b}
Dezert J.,  \emph{Combination of paradoxical sources of information within the Neutrosophic framework}, Proc. of 1st Int. Conf. on Neutrosophic Logic, Set, Prob. and Statistics, Univ. of  New Mexico, Gallup Campus, pp.~22--46, 1-3 Dec. 2001.

\bibitem{Dezert_2002}
Dezert J., \emph{An introduction to the theory of plausible and 
paradoxical reasoning}, Proc. of NM\&A 02 Conf., Borovetz, Bulgaria, Aug. 20-24, 2002. 

\bibitem{Dezert_2002b}
Dezert J., \emph{Foundations for a new theory of plausible and paradoxical 
reasoning}, Information \& 
Security, An int. Journal, edited by Prof. Tzv. Semerdjiev, CLPP, Bulgarian Acad. of 
Sci., Vol. 9, 2002.

\bibitem{Dezert_2003}
Dezert J.,  \emph{Fondations pour une nouvelle théorie du raisonnement plausible et paradoxal}, ONERA Tech. Rep. RT 1/06769/DTIM, Jan. 2003.

\bibitem{Dezert_2003d}
Dezert J., Smarandache F., \emph{On the generation of hyper-powersets for the DSmT}, Proc. of the 6th Int. Conf. on inf. fusion (Fusion 2003), Cairns, Australia, july 8-11, 2003.

\bibitem{Dezert_Smarandache_Daniel_2003}
Dezert J., Smarandache F.,  Daniel M., \emph{On the generalized pignistic transformation based on DSmT framework}, (in preparation), 2003.

\bibitem{Kennes_1991}
Kennes R., Smets Ph., \emph{Fast algorithms for Dempster-Shafer theory}, in 
Uncertainty in Knowledge Bases, B. Bouchon-Meunier, R.R. Yager, 
L.A. Zadeh (Editors), Lecture Notes in Computer Science 521, Springer-Verlag, Berlin, pp. 14-23, 1991.

\bibitem{Kennes_1992}
Kennes R., \emph{Computational Aspects of the Möbius Transformation of 
Graphs}, IEEE Trans. on SMC. 22, pp. 201-223, 1992.

\bibitem{Lefevre_2002}
Lefevre E., Colot O., Vannoorenberghe P. ``Belief functions combination and conflict management", 
Information Fusion Journal, Elsevier, 2002.

\bibitem{Shafer_1976}
Shafer G., \emph{A Mathematical Theory of Evidence}, Princeton Univ. Press, Princeton, NJ, 1976.

\bibitem{Smarandache_2002}
Smarandache F. (Editor), \emph{Proceedings of the First International Conference on 
Neutrosophics}, Univ. of New Mexico, Gallup Campus, NM, USA, 1-3 Dec. 2001, Xiquan, Phoenix, 2002.

\bibitem{Smets_1994}
Smets Ph., Kennes R., \emph{The transferable belief model}, 
Artificial Intelligence, 66(2), pp. 191-234, 1994.

\bibitem{Smets_2002}
Smets Ph., \emph{Matrix Calculus for Belief Functions}, (submitted in 2002)
{\scriptsize{\verb+http://iridia.ulb.ac.be/~psmets/MatrixRepresentation.pdf+}}.

\bibitem{Tchamova_2003}
Tchamova A., Semerdjiev T., Dezert J., \emph{Estimation of Target behavior tendencies using Dezert-Smarandache theory}, Proc. of the 6th Int. Conf. on Inf. Fusion (Fusion 2003), Cairns, Australia, july 8-11, 2003.
\end{thebibliography}
\end{document}